\documentclass[9pt,twocolumn,final]{IEEEtran}
  \usepackage{cite}
\usepackage{amsmath,amssymb,amsfonts}
\usepackage{algorithmic}
\usepackage{graphicx}
\usepackage{textcomp}
\def\BibTeX{{\rm B\kern-.05em{\sc i\kern-.025em b}\kern-.08em
    T\kern-.1667em\lower.7ex\hbox{E}\kern-.125emX}}
\numberwithin{equation}{section}

\newtheorem{prop}{Proposition}
\newtheorem{rem}{Remark}

\begin{document}

\title{A General Class of Control Lyapunov Functions and Sampled-Data Stabilization}
\author{Katerina Chrysafi and John Tsinias

\thanks{The authors are with Department of Mathematics, National Technical University of Athens, Zografou Campus 15780, Athens, Greece, e-mail: jtsin@central.ntua.gr (corresponding author), echrysafi@central.ntua.gr}}
\maketitle
\begin{abstract}
Several important contributions towards sampled-data and hybrid feedback stabilization have appeared in the literature. The present work extends recent results by second author concerning sampled-data feedback stabilization for affine in the control of nonlinear systems with nonzero drift term, under the presence of a generalized control Lyapunov function associated with appropriate Lie algebraic hypotheses concerning the dynamics of the system. The main results of present work, constitute a generalization of the well-known "Artstein-Sontag" theorem on  asymptotic stabilization by means of an almost smooth feedback controller. The analysis is limited to the affine single-input nonlinear systems with nonzero drift term, however, the results can easily be extended to the multi-input case. An illustrative example of the derived results is included.
\end{abstract}
\section{Introduction}
 Within the framework of problems concerning the stabilization of nonlinear autonomous systems by means of sampled-data and hybrid feedback control, several important results have appeared in the literature; see for instance \cite{art:1}, \cite{art:3} -\cite{art:7}, \cite{art:9} -\cite{art:18} and relative references therein. In the works \cite{art:21} - \cite{art:24}, various concepts of sampled-data stabilization are introduced for general systems
\begin{equation}
\dot{x}=f(x,u),\,\, (x,u)\in \mathbb{R}^{\text{n}}\times \mathbb{R}^{\text{m}}\,,\,
f(0,0)=0
\label{1.1}
\end{equation}
and Lyapunov type sufficient conditions are established guaranteeing Sampled-Data Feedback Semi Global Asymptotic Stabilization (SDF-SGAS). Especially, for the case of affine in the control systems:
\begin{equation} {\dot{x}=F(x,u):=f(x)+ug(x),\,(x,u)\in\mathbb{R}^{\text{n}} \times \mathbb{R}} \,,\, {f\left(0\right)=0} \label{1.2}
\end{equation}
Lie algebraic sufficient conditions are derived in \cite{art:21}, \cite{art:23} and \cite{art:24} for SDF-SGAS. These conditions constitute extensions of the familiar ``Artstein-Sontag'' Lyapunov-like sufficient condition for asymptotic stabilization of systems \eqref{1.2} by means of an almost smooth feedback. (see \cite{art:2} , \cite{art:19} and \cite{art:20}).
\par In the present paper we extend the main result of \cite{art:24} by deriving a general set of sufficient Lie algebraic conditions for SDF-SGAS for systems \eqref{1.2} and the additional possibility of stabilizing \eqref{1.2} by means of a \emph{bounded} sampled-data feedback. We present below the concept of SDF-SGAS as introduced in \cite{art:24}. We assume that $f:{\mathbb R}^{\text{n}} \times {\mathbb R}^{\text{m}} \to {\mathbb R}^{\text{n}} $ is Lipschitz continuous and we denote by $\pi(\cdot )=\pi(\cdot ,s,x,u)$ the trajectory of \eqref{1.1} with initial condition $\pi(s,s,x ,u)=x \in {\mathbb R}^{\text{n}} $ corresponding to certain measurable and locally essentially bounded control $u:[s,T_{\max } )\to {\mathbb R}^{\text{m}} $, where $T_{\max } =T_{\max } (s,x ,u)$ is the corresponding maximal existing time of the trajectory. We say that system \eqref{1.1} is \emph{Semi-Globally Asymptotically Stabilizable by Sampled-Data Feedback  (SDF-SGAS)}, if for every $R>0$ and for any given partition of times $T_{1} :=0<T_{2} <T_{3} <... <T_{\nu } <... {\kern 1pt} {\kern 1pt} \, \, ,{\rm with}\, \, \,  T_{\nu } \to \infty$ and $T_{\nu+1}-T_{\nu}\,,\,\nu=1,2,...$ bounded, there exist a neighborhood $\Omega $ of zero with $B[0,R]:=\left\{x\in {\mathbb R}^{\text{n}} :|x|\le R\right\}\subset \Omega $; ($|x|$ denotes the usual Euclidean norm of the vector $x\in\mathbb{R}^{\text{n}}$) and a map $k:{\mathbb R}^{+} \times \Omega \to  {\mathbb R}^{\text{m}} $ such that for any $x\in \Omega $ the map $k(\cdot ,x): {\mathbb R}^{+} \to  {\mathbb R}^{\text{m}} $ is measurable and essentially bounded and the trajectory $\pi(\cdot )$ of the sampled-data closed loop system $\dot{\pi}=f(\pi,k(t,\pi(T_{i} )))$, for $t\in [T_{i} ,\; T_{i+1} ),\; \, \, i=1,2,...$ with initial  $\pi(0)\in \Omega $ satisfies the following properties:
 Stability: For every $\varepsilon >0$ there exists $\delta =\delta (\varepsilon )>0$ such that $|\pi(0)|\leq \delta\Rightarrow |\pi(t)|\leq\varepsilon$, for any $t\geq0$ and for any  $\pi(0)\in \Omega$; Attractivity: 
$\mathop{\lim }\limits_{t\to +\infty } \pi(t)=0$, for any initial $\pi(0)\in \Omega$. We say that system \eqref{1.2} is \emph{bounded SDF-SGAS (BSDF-SGAS)}, if it is SDF-SGAS and in addition there exists a constant $C=C_{\Omega}>0$ such that the corresponding map $k:\mathbb{R}^{+}\times\Omega\to\mathbb{R}^{\text{m}}$ satisfies : $|k(t,x)|\leq C$, $\forall x\in\Omega$, $t\geq0$ near zero. The following result constitutes a direct extension of [24, Proposition 2].
\begin{prop}
For the system \eqref{1.1} assume that there exists a positive definite $C^{0}$ function $V:\mathbb{R}^{\normalfont \text{n}} \to\mathbb{R}^{+}$ and a function $a\in K$   (namely, $a(\cdot)$  is continuous and strictly increasing with $a(0)=0$ ) such that for every $\sigma>0$ and $x\neq0$  there exist a constant $\varepsilon=\varepsilon(x)\in(0,\sigma]$   and a measurable and locally essentially bounded control $u(\cdot,x):[0,\varepsilon]\to\mathbb{R}^{ \normalfont \text{m}} $  satisfying
\begin{subequations}
\label{1.3}
\begin{equation}
\label{1.3a}
 V(\pi(\varepsilon ,0,x ,u(\cdot ,x) ))<V(x);  
\end{equation}
\begin{equation}
\label{1.3b}
 V(\pi(s,0,x ,u(\cdot ,x) ))\le a(V(x)),\; \; \forall s\in [0,\varepsilon ]  
\end{equation}
\end{subequations}
Then, system \eqref{1.1} is SDF-SGAS. If, in addition to \eqref{1.3a} and \eqref{1.3b} we assume that  for any bounded   nonempty  neighborhood $\Omega$ of zero $0\in\mathbb{R}^{\normalfont \text{n}}$ there exists a constant $C=C_{\Omega}>0$ such that the corresponding value $u(\cdot, \cdot)$  satisfies
\begin{equation}
    \label{1.4}
    |u(t,x)|\leq C,\,\,\forall x\in\Omega,\,\,t\geq0\,\,\text{near zero}
\end{equation}
then system \eqref{1.1} is  BSDF-SGAS.
\end{prop}
\par The paper is organized as follows : Section II contains the precise statements of main results (Proposition 2 and 3), Section III contains the proofs of our results and in Section IV an illustrative example is examined.
\par We require some standard notations and definitions. For any pair of $C^{1}$ mappings $X:\mathbb{R}^{\text{n}}\to\mathbb{R}^{\text{k}}$, $Y:\mathbb{R}^{\text{k}}\to\mathbb{R}^{\text{l}}$ we adopt the notation $XY:=(DY)X$, with $DY$ being the derivative of $Y$. By $[\cdot,\cdot]$ we denote the Lie bracket operator, namely, $[X,Y]=XY-YX$ for any pair of $C^{1}$ mappings $X,Y:\mathbb{R}^{\text{n}}\to\mathbb{R}^{\text{n}}$. By $\text{Lie}\left\{f,g\right\}$ we denote the Lie Algebra generated by $\left\{f,g\right\}$. We also require some elementary concepts concerning the order of a vector field $\Delta\in\text{Lie}\left\{f,g\right\}$: Let $L_{1}:=\text{span}\left\{f,g\right\}$ and $L_{i+1}:=\text{span}\left\{[X,Y],\,X\in L_{i},\, Y\in L_{1}\right\}$, $i=1,2,...$. Then, for any nonzero $\Delta\in\text{Lie}\left\{f,g\right\}$ define:
$\text{order}_{\left\{f,g\right\}}\Delta:=1 \, \text{ if }\, \Delta\in L_{1}\setminus \left \{ 0 \right \}$ and $\text{order}_{\left\{f,g\right\}}\Delta:=k>1\, \text{ if }\, \Delta=\Delta_{1}+\Delta_{2}\,
 \text{with}\,\,\Delta_{1}\in L_{k}\setminus\left \{ 0 \right \}  
\text{and}\,\,\Delta_{2}\in\text{span}\left \{\bigcup_{i=1}^{i=k-1}L_{i}  \right \}$.
For simplicity, we adopt the notation $\text{order}\Delta=\text{order}_{\left\{f,g\right\}}\Delta$. Finally, we use the notation $\text{order}_{\left\{g\right\}}\Delta=l$, $l\in\mathbb{N}_{0}$ for the case where $\Delta\in\text{Lie}\left\{f,g\right\}$ is contained to the linear span of all Lie monomials of $f$ and $g$ involving $g$ exactly $l$ times.
\section{Main Results}
We assume that the dynamics  of system \eqref{1.2} are smooth  ($C^{\infty}$)(although, our main results can directly be extended for systems under weaker regularity assumptions). In order to present our new results we first introduce a subalgebra of $\text{Lie}\left\{f,g\right\}$, which plays the central role in order to provide our Lie-algebraic conditions for stabilization of \eqref{1.2} by means of a sampled- data feedback. Consider the following vector fields: 
\begin{equation}
\label{2.1}
\begin{split}
    &\lambda_{1,0}:=f
    \\
    &\lambda_{\kappa,j}:=
    \\
    &\sum_{\substack{r_{1},...,r_{j}\in \mathbb{N}_{0}, \\ r_{1}+...+r_{j}=\kappa-j-1}}\begin{Bmatrix}
[...[[[...[[...[[...[[f,g],
\underset{r_{j}}{\underbrace{f],...,f}}],g],...\\
...,g], \underset{r_{2}}{\underbrace{f],...,f}}],g], \underset{r_{1}}{\underbrace{f],...,f}}]   
\end{Bmatrix}
\\
&\kappa=2,3,...;\,\,j=1,...,\kappa-1
\end{split}
\end{equation}
For instance, according to definition \eqref{2.1}: $\lambda_{2,1}=[f,g]$, $\lambda_{3,1}=[[f,g],f]$, $\lambda_{3,2}=[[f,g],g]$, $\lambda_{4,1}=[[[f,g],f],f]$, $\lambda_{4,2}=[[[f,g],f],g]+[[[f,g],g],f]$, $\lambda_{4,3}=[[[f,g],g],g]$, $\lambda_{5,1}=[[[[f,g],f],f],f]$, $\lambda_{5,2}=[[[[f,g],f],f],g]+[[[[f,g],f],g],f]+[[[[f,g],g],f],f]$, $\lambda_{5,3}=[[[[f,g],f],g],g]+[[[[f,g],g],f],g]+[[[[f,g],f],f],g]$, $\lambda_{5,4}=[[[[f,g],f],f],f]$ etc.
\\
\par Define the following subalgebra of $\text{Lie}\left\{f,g\right\}$:
\begin{equation}
    \label{2.2}
    L\left\{f,g\right\}:=\text{span}\left\{\lambda_{\kappa,j},\,\,\kappa=1,2,...,\,j=0,...,\kappa-1\right\}
\end{equation}
Obviously, $L\left\{f,g\right\}\subset\text{Lie}\left\{f,g\right\}$ and, according to \eqref{2.1} and the notations of previous section, we have $\text{order}\lambda_{\kappa,j}=\kappa\geq1$ and $\text{order}_{\left\{g\right\}}\lambda_{\kappa,j}=j\geq0$. The following result provides Lie algebraic conditions for SDF-SGAS for systems \eqref{1.2} in terms of the elements of $L\left\{f,g\right\}$ and generalizes Proposition 3 in \cite{art:24}.
\begin{prop}
For the system \eqref{1.2} assume that there exists a smooth function $V:\mathbb{R}^{\normalfont \text{n}}\to\mathbb{R}^{+}$ being positive definite and proper, such that for every $x\neq0$, either $(gV)(x)\neq0$, or one of the following properties is fulfilled: Either
\begin{equation}
    \label{2.4}
    (gV)(x)=0\Rightarrow(fV)(x)<0
\end{equation}
or there exists an integer $N=N(x)\geq1$ such that 
\begin{equation}
    \label{2.5}
    \begin{split}
    &(gV)(x)=0
    \\
    &\Rightarrow\begin{Bmatrix}
(\Delta_{1}\Delta_{2}... \Delta_{k} V)(x)=0,\,\,k=1,2,... \\
\forall \Delta _{1} , \Delta_{2}... , \Delta_{k} \in L\{ f,g\}:\,\sum _{p=1}^{k}\mathrm{order}\Delta_{p}  \le N \\ 
\end{Bmatrix}
\end{split}
\end{equation}
and in such a way that one of the following two properties is fulfilled:
\begin{equation}\label{2.6}\text{Property P1}:\,\,\,\,\,\,(f^{N+1} V)(x)<0\end{equation}
\begin{equation}
    \label{2.7}
\text{Property P2}:\,\,\,\,\,\,(f^{N+1} V)(x)\leq0    
\end{equation}
where $(f^{i}V)(x):=f(f^{i-1}V)(x)$, $i=2,3,...$, $(f^{1}V)(x):=(fV)(x)$ and further, one of the following holds:
\\
\textbf{P2(i)}: $N$ is a positive integer and there exists an \textbf{odd} integer $j=j(N)\in\{1,...,N\}$ such that
\begin{equation}
    \label{2.8}
    (\lambda_{N+1,j}V)(x)\neq0
\end{equation}
where $\lambda_{N+1,i}$, $i\in\{1,...,N\}$ are generators of the Lie subalgebra $L\{f,g\}$ defined by \eqref{2.1} and \eqref{2.2}, respectively. Moreover, in addition to \eqref{2.5} assume that:
\begin{subequations}
\label{2.9}
\begin{equation}
    \label{2.9a}
    \begin{split}
     &(gV)(x)=0\\
     &\Rightarrow\begin{Bmatrix}
(\Delta_{1}\Delta_{2}... \Delta_{k} V)(x)=0,\,\,k=1,2,... \\
\forall \Delta _{1} , \Delta_{2}... , \Delta_{k} \in L\{ f,g\}:
\\
\sum _{p=1}^{k}\mathrm{order}\Delta_{p}=N+1;\,\sum_{p=1}^{k}\mathrm{order}_{\{g\}}\Delta_{p}=q \\ 
\end{Bmatrix}
\end{split}
\end{equation}
\begin{equation}
    \label{2.9b}
    \text{for every \textbf{even} positive integer}\,q<j,\text{provided that}\,j>2.
\end{equation}
\end{subequations}
\textbf{P2(ii)}: $N$ is an \textbf{odd} integer $N=N(x)>2$ and there exists an \textbf{odd} integer $j=j(N)\in\{1,...,N-2\}$ such that \eqref{2.8} holds and implication \eqref{2.9a} is fulfilled for every \textbf{even} positive integer $q:j<q<N$.
\\
\textbf{P2(iii)} $N$ is \textbf{even} and 
\begin{equation}
    \label{2.10}
    (\lambda_{N+1,N}V)(x)<0
\end{equation}
Then, system \eqref{1.2} satisfies conditions \eqref{1.3a} and \eqref{1.3b} of Proposition 1, hence, it is SDF-SGAS.
\end{prop}
\begin{rem} The conditions imposed in Proposition 2 are weaker than the corresponding assumptions in [24, Proposition 3]. One difference is that in all assumptions made in the statement of Proposition 2 the Lie subalgebra $L\{f,g\}$ is involved, instead of $\text{Lie}\{f,g\}\setminus\{g\}$  considered in \cite{art:24}. Another essential difference is the validity of assumption \eqref{2.8}, instead of the stronger condition $(\lambda_{N+1,N}V)(x)\neq0$, where $N=N(x)\geq1$, imposed in \cite{art:24}.
\end{rem}
\begin{rem} Condition \eqref{2.9} holds, if for instance, the following is fulfilled: For every $x\neq0$ with $(gV)(x)=0$ and positive even $q<j$ it holds: $\mathrm{span}\{(\Delta_{1}\Delta_{2}... \Delta_{k} V)(x),\,\,k=1,2,...:\Delta _{1} , \Delta_{2}... , \Delta_{k} \in L\{ f,g\}$, with $\sum_{p=1}^{k}\mathrm{order}\Delta_{p}=N+1;\,\,\,\,\sum_{p=1}^{k}\mathrm{order}_{\{g\}}\Delta_{p}=q\}=\mathrm{span}\{(\Delta_{1}\Delta_{2}... \Delta_{k} V)(x),\,\,k=1,2,...:
\Delta _{1} , \Delta_{2}... , \Delta_{k} \in L\{ f,g\}$, with $\sum _{p=1}^{k}\mathrm{order}\Delta_{p}  \leq N,\,\, 
\sum_{p=1}^{k}\mathrm{order}_{\{g\}}\Delta_{p}=q-1\}$.
Indeed, the previous equality in conjunction with \eqref{2.5} implies \eqref{2.9}. It  is worthwhile noticing that the assumption above, is a variant of the Hermes condition (see \cite{art:8}), namely, that when $q>0$ is even, it is assumed that for $x\neq0$, the span of all  Lie monomials of $f$ and $g$  involving $g$ at most $q-1$ times evaluated at $x$ is equal to the span of the corresponding  Lie monomials of $f$ and $g$ involving $g$ at most $q$ times.
  \end{rem}
\par The following result provides Lie algebraic conditions for BSDF-SGAS for systems \eqref{1.2} in terms of the subalgebra $L\{f,g\}$.
\begin{prop}
For the system \eqref{1.2} assume that there exists a smooth function $V:\mathbb{R}^{\normalfont \text{n}}\to\mathbb{R}^{+}$  being positive definite and proper such that all conditions of Proposition 2 are fulfilled and further we assume that 
\begin{itemize}
\item in addition to \eqref{2.4}, there exists a pair of continuous nonnegative functions $\theta:\mathbb{R}^{\normalfont \text{n}}\to\mathbb{R}^{+}$ and  $\xi:\mathbb{R}^{+}\to\mathbb{R}^{+}$ such that
\begin{equation}
    \label{2.11}
    |(fV)(\omega)+\theta(\omega)|\leq\xi(|\omega|)|(gV)(\omega)|,\,\,\forall \omega\in\mathbb{R}^{\normalfont \text{n}}
\end{equation}
 \item Property P2(ii) is strengthened by assuming that is fulfilled with $j=N$ (odd) and in addition that the following implication is imposed for the specific odd $N$:
\begin{equation}
    \label{2.12}
    \begin{split}
     &(gV)(x)=0\\
     &\Rightarrow\begin{Bmatrix}
(\Delta_{1}\Delta_{2}... \Delta_{k} V)(x)=0,\,\,k=1,2,... \\
\forall \Delta _{1} , \Delta_{2}... , \Delta_{k} \in L\{ f,g\}:
\\
\sum _{p=1}^{k}\mathrm{order}\Delta_{p}=N+1;
\\
\sum_{p=1}^{k}\mathrm{order}_{\{g\}}\Delta_{p}\leq N-1 
\end{Bmatrix}
\end{split}
\end{equation}
 \item Property P2(iii) is strengthened by assuming that  the corresponding even integer $N$ satisfies the additional implication \eqref{2.12}.
 \end{itemize}
Then, system \eqref{1.2} is BSDF-SGAS.
\end{prop}
\section{Proof of main results}
\textbf{Proof of Proposition 2:} We first select a continuous nonnegative function $\theta:\mathbb{R}^{\text{n}}\to\mathbb{R}^{+}$ and define:
\begin{equation}
    \label{3.1}
    u=u(x):=\begin{cases}
-\left ( \frac{(fV)(x)+\theta(x)}{((gV)(x))^{2}}+1 \right )(gV)(x),\, (gV)(x)\neq0 \\ 
 0,\,\,\,\,\,\,(gV)(x)=0\,\text{and}\,(fV)(x)<0 
\end{cases}
\end{equation}
Then, it can easily be established that \eqref{3.1} together with \eqref{2.4} imply that for every sufficiently small $\varepsilon>0$  the corresponding trajectory of \eqref{1.2}, with $u$ as defined by \eqref{3.1}, satisfies:
\begin{subequations}
\label{3.2}
\begin{equation}
    \label{3.2a}
    V(\pi(s,0,x,u))<V(x),\,\,\forall s\in(0,\varepsilon]
\end{equation}
\begin{equation}
    \label{3.2b}
    \begin{split}
    &\forall x\neq0:\text{such that either}\,\,(gV)(x)\neq0\\
    &\mathrm{or}\,\,(gV)(x)=0\,\,\mathrm{and}\,\,(fV)(x)<0
    \end{split}
\end{equation}
\end{subequations}
Assume next that there exists an integer $N=N(x)\geq1$ satisfying \eqref{2.5} together with one of the Properties P1, P2. We apply an extension of the procedure employed for the proof of [24, Proposition 3]. Define:
\begin{equation}
    \label{3.3}
    X:=f+u_{1}g,\,\,Y:=f+u_{2}g
\end{equation}
and denote by $X_{t}(z)$ and $Y_{t}(z)$ the trajectories of the systems $\dot{x}=X(x)$ and $\dot{y}=Y(y)$ respectively, initiated at time $t=0$ from some $z\in\mathbb{R}^{\text{n}}$. For every constant $\rho>0$ define:
\begin{equation}
    \label{3.4}
    R(t):=\left ( X_{\rho t}\circ Y_{t} \right )(x),\,\,t\geq0,\,\,R(0)=x
\end{equation}
\begin{equation}
    \label{3.5}
    m(t):=V(R(t)),\,\,t\geq0
\end{equation}
and denote by $m^{(n)}(\cdot)$, $n=1,2,...$ its $n$-time derivative. We prove that for every nonzero $x$, condition \eqref{2.5}, together with one of the rest properties imposed in the statement of Proposition 2, imply existence of a constant $\rho=\rho(x)>0$ and a pair of constant inputs $u_{1}=u_{1}(x)$ and $u_{2}=u_{2}(x)$ such that
\begin{equation}
\label{3.6}
    m^{(n)}(0)=0,\,\,n=1,...,N
\end{equation}
\begin{equation}
    \label{3.7}
    m^{(N+1)}(0)<0
\end{equation}
In order to establish \eqref{3.6} and \eqref{3.7}, we express the time derivatives $m^{(n)}(\cdot)$, $n=1,2,...$ of the map $m(\cdot)$ defined by \eqref{3.5} in terms of the elements of the Lie subalgebra $L\{f,g\}$ defined by \eqref{2.1}, \eqref{2.2}. Indeed, we take into account definitions \eqref{3.3}-\eqref{3.5} and apply the Campbell-Baker-Hausdorff formula to the right hand side of \eqref{3.4}. Then, as in the  proof of [24, Proposition 3], by setting
\begin{equation}
\label{3.8}
    u_{2}:=-\rho u_{1},\,\,\rho>0
\end{equation}                                                      
we find:
\begin{subequations}
\label{3.9}
\begin{equation}
    \label{3.9a}
    m^{(1)}(0)=(\rho+1)(fV)(x),\,\,\mathrm{for}\,\,n=1
\end{equation}
\begin{equation}
    \label{3.9b}
    \begin{split}
    m^{(2)}(0)=&(\rho+1)^{2}(f^{2}V)(x)+u_{1}\rho(\rho+1)(\lambda_{2,1}V)(x),\\
    &\mathrm{for}\,\,n=2
    \end{split}
\end{equation}
\begin{equation}
\label{3.9c}
\begin{split}
 &m^{(n)}(0)=(\rho+1)^{n}(f^{n}V)(x)
 \\
 &+\sum_{i=1}^{n-2}u_{1}^{i}\left ( \Pi_{n,i}(\rho;x)+\rho^{n-1}(\rho+1)(\lambda_{n,i}V)(x) \right )
 \\
 &+u_{1}^{n-1}\rho^{n-1}(\rho+1)(\lambda_{n,n-1}V)(x),\,\,\mathrm{for}\,\,n=3,4,...
 \end{split}
\end{equation}
\end{subequations}
where each map $\Pi_{n,i}(\rho;x)$ for $n\geq3$ and $i=1,...,n-2$, appearing on the right hand side of \eqref{3.9c}, satisfies the following properties:
\\
I.It is exclusively dependent on $\rho\in\mathbb{R}$ and $x\in\mathbb{R}^{\text{n}}$, thus, it is independent of $u$.
\\
II.For each fixed  $x\in\mathbb{R}^{\text{n}}$ the map $\Pi_{n,i}(\rho;x)$ is a nonzero polynomial with respect to $\rho$ in such a way that
\begin{subequations}
\label{3.10}
\begin{equation}
\label{3.10a}
\mathrm{deg}\Pi_{n,i}(\cdot;x)=n,\,\,n\geq3,\,\,i=1,...,n-2    
\end{equation}
\begin{equation}
    \label{3.10b}
    \begin{split}
    \Pi_{n,i}(\rho;x)\,\,&\text{and}\,\,q(\rho)=\rho^{n-1}(\rho+1),\,\,n\geq3\\
    &\text{are linearly independent}
    \end{split}
\end{equation}
\end{subequations}
III.For every $\rho\in\mathbb{R}$ and $x\in\mathbb{R}^{\text{n}}$ the following holds: 
\begin{equation}
\label{3.11}
\begin{split}
    &\Pi_{n,i}(\rho;x)
    \\
    &\in\mathrm{span}\begin{Bmatrix}
(\Delta_{1}\Delta_{2}... \Delta_{k} V)(x),\,\,k=1,2,..., \\
\Delta _{1} , \Delta_{2}... , \Delta_{k} \in L\{ f,g\}:
\\
\sum _{p=1}^{k}\mathrm{order}\Delta_{p}=n;\, \sum_{p=1}^{k}\mathrm{order}_{\{g\}}\Delta_{p}=i \\
\end{Bmatrix},\\
&\forall i,n\in\mathbb{N}:n\geq3,\,\,1\leq i\leq n-2
\end{split}
\end{equation}
By taking into account definitions \eqref{2.1}, \eqref{2.2}, assumption \eqref{2.5} and applying \eqref{3.11} with $n=N$, we find:
\begin{equation}
    \label{3.12}
    \Pi_{N,i}(\rho;x)=0,\,\,i=1,...,N-2\,\,\text{for all}\,\,\rho\in\mathbb{R}
\end{equation}
\begin{equation}
    \label{3.13}
    \begin{split}
    (f^{i}V)(x)&=0,\,\,i=1,...,N;\\
    (\lambda_{N,j}V)(x)&=0,\,\,j=1,...,N-1
    \end{split}
\end{equation}
From \eqref{3.12}, \eqref{3.13} and validity of \eqref{3.9a}, \eqref{3.9b}, \eqref{3.9c} with $n=N$, it follows that \eqref{3.6} holds. Next, we establish \eqref{3.7}. We distinguish four cases:
\\
 \textbf{CASE 1}: Suppose that Property P1 is satisfied, namely, condition \eqref{2.6} holds. Then, by recalling \eqref{3.9b} and \eqref{3.9c} with $n=N+1$ and by setting $u_{1}=0$, $u_{2}=-\rho u_{1}=0$, it follows that the desired \eqref{3.7} holds for every choice of $\rho>0$.
 \\
 \textbf{CASE 2}: Suppose next that both \eqref{2.7} and P2(i) hold. Consider first the case $N=1$. We recall assumption \eqref{2.8} with $j=1$, namely:
 \begin{equation}
     \label{3.14}
     (\lambda_{2,1}V)(x)\neq0
 \end{equation}
From  \eqref{3.9b}, \eqref{3.14} and validity of \eqref{2.7}, it follows that for every $\rho>0$ there exists an arbitrarily small constant $u_{1}=u_{1}(x)\neq0$ such that inequality \eqref{3.7} holds with $N=1$.
\\
Next, we establish \eqref{3.7} for the case $N\geq2$, under \eqref{2.7} and \eqref{2.8}. By invoking \eqref{3.9c} with $n=N+1$ we get:
\begin{equation}
\label{3.15}
 \begin{split}
 &m^{(N+1)}(0)=(\rho+1)^{N+1}(f^{N+1}V)(x)
 \\
 &+\sum_{i=1}^{N-1}u_{1}^{i}\left ( \Pi_{N+1,i}(\rho;x)+\rho^{N}(\rho+1)(\lambda_{N+1,i}V)(x) \right )
 \\
 &+u_{1}^{N}\rho^{N}(\rho+1)(\lambda_{N+1,N}V)(x)
 \end{split}   
\end{equation}
By taking into account \eqref{2.9a}, \eqref{2.9b} and \eqref{3.11} with $n=N+1$ and $i=1,...,N-1$, it follows that each mapping $\Pi_{N+1,i}(\rho;x)$, $N\geq2$, $i=1,...,N-1$ above satisfies:
\begin{equation}
    \label{3.16}
    \begin{split}
    &\Pi_{N+1,q}(\rho;x)=0,\,\text{for every even}
    \\
    &\text{positive integer}\,q<j\,\text{and for every}\,\rho\in\mathbb{R}
    \end{split}
\end{equation}
where $j\in\{1,...,N\}$ is the odd integer satisfying \eqref{2.8}. Next, we recall definitions \eqref{2.1} and \eqref{2.2}, which imply:
\begin{equation}
    \label{3.17}
    \begin{split}
&(\lambda_{n,i}V)(x)
\\
&\in\mathrm{span}\begin{Bmatrix}
(\Delta_{1}\Delta_{2}... \Delta_{k} V)(x),\,\,k=1,2,..., \\
\Delta _{1} , \Delta_{2}... , \Delta_{k} \in L\{ f,g\}:
\\
\sum _{p=1}^{k}\mathrm{order}\Delta_{p}=n;\, \sum_{p=1}^{k}\mathrm{order}_{\{g\}}\Delta_{p}=i 
\end{Bmatrix},\\
&\text{for}\,\,n\geq2,\,i=1,...,n-1
\end{split}
\end{equation}
Then, by taking into account \eqref{3.17} with $n=N+1$ and recalling assumption \eqref{2.9a}, \eqref{2.9b} it follows:
\begin{equation}
    \label{3.18}
    (\lambda_{N+1,q}V)(x)=0,\,\,\text{for every even}\,q<j
\end{equation}
From \eqref{3.15}, \eqref{3.16} and \eqref{3.18} we obtain:
\begin{equation}
    \label{3.19}
 \begin{split}
 &m^{(N+1)}(0)=(\rho+1)^{N+1}(f^{N+1}V)(x)
 \\
 &+\sum_{k=0}^{\frac{j-1}{2}}u_{1}^{2k+1}\Big( \Pi_{N+1,2k+1}(\rho;x)+\rho^{N}(\rho+1)(\lambda_{N+1,2k+1}V)(x) \Big)
 \\
 &+\sum_{k=j+1}^{N-1}u_{1}^{k}\left ( \Pi_{N+1,k}(\rho;x)+\rho^{N}(\rho+1)(\lambda_{N+1,k}V)(x) \right )
 \\
 &+u_{1}^{N}\rho^{N}(\rho+1)(\lambda_{N+1,N}V)(x),\,\mathrm{for}\,\,1\leq j\leq N
 \end{split}
\end{equation}
We are in a position to establish \eqref{3.7} by induction as follows: 
\\
\textbf{Step 1}: Since, according to \eqref{3.10b} with $n=N+1$ and $i=1$, the polynomials $\Pi_{N+1,1}(\rho;x)$ and $q(\rho)=\rho^{N}(\rho+1)$ are linearly independent, we examine two subcases:
\\
\textbf{\emph{Subcase i}}: There exists a constant $\rho=\rho(x)\in(0,1]$ with $\Pi_{N+1,1}(\rho;x)+\rho^{N}(\rho+1)(\lambda_{N+1,1}V)(x)\neq0$. 
Note that, since $j$ is odd, each $2k+1$, $k=0,1,...,(j-1)/2$ is an odd positive integer. It then follows from \eqref{3.19} that there exist constants $u_{1}=u_{1}(x)$, $u_{2}=-\rho u_{1}$, with $|u_{1}|, |u_{2}|$ arbitrarily small, satisfying \eqref{3.7}.
\\
\textbf{\emph{Subcase ii}}: $\Pi_{N+1,1}(\rho;x)=0$, $\forall \rho\in\mathbb{R}$; $(\lambda_{N+1,1}V)(x)=0$. Then \eqref{3.19} is rewritten:
\begin{equation}
    \label{3.20}
    \begin{split}
 &m^{(N+1)}(0)=(\rho+1)^{N+1}(f^{N+1}V)(x)
 \\
 &+\sum_{k=1}^{\frac{j-1}{2}}u_{1}^{2k+1}\Big( \Pi_{N+1,2k+1}(\rho;x)+\rho^{N}(\rho+1)(\lambda_{N+1,2k+1}V)(x) \Big)
 \\
 &+\sum_{k=j+1}^{N-1}u_{1}^{k}\left ( \Pi_{N+1,k}(\rho;x)+\rho^{N}(\rho+1)(\lambda_{N+1,k}V)(x) \right )
 \\
 &+u_{1}^{N}\rho^{N}(\rho+1)(\lambda_{N+1,N}V)(x),\,\text{provided that}\,3\leq j\leq N
 \end{split}
\end{equation}
 and we proceed with the next step:
 \\
\textbf{Step 2}: Since, according to \eqref{3.10b} with $n=N+1$ and $i=3$, the polynomials $\Pi_{N+1,3}$ and $q(\rho)=\rho^{N}(\rho+1)$ are linearly independent, we again distinguish two subcases :
\\
\textbf{\emph{Subcase i}}: There exists a constant $\rho=\rho(x)\in(0,1]$ such that $\Pi_{N+1,3}(\rho;x)+\rho^{N}(\rho+1)(\lambda_{N+1,3}V)(x)\neq0$. It then follows from the latter fact and \eqref{3.20} that there exist constants $u_{1}=u_{1}(x)$, $u_{2}=-\rho u_{1}$, with $|u_{1}|, |u_{2}|$ arbitrarily small, satisfying \eqref{3.7}. 
\\
\textbf{\emph{Subcase ii}}: $\Pi_{N+1,3}(\rho;x)=0$, $\forall \rho\in\mathbb{R}$; $(\lambda_{N+1,3}V)(x)=0$. Then \eqref{3.20} becomes:
\begin{equation}
    \label{3.21}
    \begin{split}
 &m^{(N+1)}(0)=(\rho+1)^{N+1}(f^{N+1}V)(x)
 \\
 &+\sum_{k=2}^{\frac{j-1}{2}}u_{1}^{2k+1}\Big( \Pi_{N+1,2k+1}(\rho;x)+\rho^{N}(\rho+1)(\lambda_{N+1,2k+1}V)(x) \Big)
 \\
 &+\sum_{k=j+1}^{N-1}u_{1}^{k}\left ( \Pi_{N+1,k}(\rho;x)+\rho^{N}(\rho+1)(\lambda_{N+1,k}V)(x) \right )
 \\
 &+u_{1}^{N}\rho^{N}(\rho+1)(\lambda_{N+1,N}V)(x),\,\text{provided that}\,\,5\leq j\leq N
 \end{split}
\end{equation}
and we proceed with the next step. Particularly, by taking into account \eqref{3.21} we may proceed quite similarly by induction and conclude that, either \eqref{3.7} holds with $|u_{1}|$, $|u_{2}|$ arbitrarily  small, satisfying \eqref{3.8} for some $\rho>0$, or
\begin{subequations}
\label{3.22}
\begin{equation}
    \label{3.22a}
    \begin{split}
    &\Pi_{N+1,i}(\rho;x)=0,\,\,\forall\,\rho\in\mathbb{R},\\
    &\text{for every odd}\,i\in\{1,3,...,j-2\}
    \end{split}
\end{equation}
\begin{equation}
    \label{3.22b}
    (\lambda_{N+1,i}V)(x)=0,\,\,\text{for every odd}\,i\in\{1,3,...,j-2\}
\end{equation}
\end{subequations}
and the procedure is terminated at the Step $(j+1)/2$  below, where, according to \eqref{3.22a} and \eqref{3.22b}, the original  expression \eqref{3.15} is written:
\begin{equation}
    \label{3.23}
\begin{split}
 &m^{(N+1)}(0)=(\rho+1)^{N+1}(f^{N+1}V)(x)
 \\
 &+\sum_{k=j}^{N-1}u_{1}^{k}\Big( \Pi_{N+1,k}(\rho;x)+\rho^{N}(\rho+1)(\lambda_{N+1,k}V)(x) \Big)
 \\
 &+u_{1}^{N}\rho^{N}(\rho+1)(\lambda_{N+1,N}V)(x)
 \end{split}
\end{equation}
\textbf{Step} \boldmath $(j+1)/2$ \unboldmath: Due to \eqref{2.8}, only one case is examined. Particularly, since $j$ satisfies \eqref{2.8}, then \eqref{3.23} and validity of \eqref{2.7}  imply that for every $\rho\in(0,1]$, there exist constants $u_{1}=u_{1}(x)$, $u_{2}=-\rho u_{1}$ with $|u_{1}|$, $|u_{2}|$ arbitrarily small, such that \eqref{3.7} holds.
\\
\textbf{CASE 3}: Assume next that Property P2(ii) is imposed. Namely, $N$ is odd, \eqref{2.7} holds and there exists an odd integer $j=j(N)\in\{1,...,N\}$ such that \eqref{2.8} and \eqref{2.9a}, for all even positive integers $q:j<q<N$,  are fulfilled. The case where $N=1$ has already been examined in CASE 2. Using \eqref{3.15}, we perform the opposite procedure with this employed in CASE 2 to establish \eqref{3.7}. By taking into account the fact that \eqref{2.9a} holds for all even positive integers $q:j<q<N$ and the fact that both $N,j$ are odd, it follows that \eqref{3.15} is rewritten as:
\begin{equation}
\label{3.24}
 \begin{split}
 &m^{(N+1)}(0)=(\rho+1)^{N+1}(f^{N+1}V)(x)
 \\
 &+\sum_{k=1}^{j}u_{1}^{k}\Big( \Pi_{N+1,k}(\rho;x)+\rho^{N}(\rho+1)(\lambda_{N+1,k}V)(x) \Big)
 \\
 &+\sum_{k=\frac{j+3}{2}}^{\frac{N-1}{2}}u_{1}^{2k-1}\Big( \Pi_{N+1,2k-1}(\rho;x)
 \\
 &+\rho^{N}(\rho+1)(\lambda_{N+1,2k-1}V)(x) \Big)+u_{1}^{N}\rho^{N}(\rho+1)(\lambda_{N+1,N}V)(x)
 \end{split}
\end{equation}
where, for simplicity, we have assumed here that $N-j\geq4$.
\\
\textbf{Step 1}: For the expression \eqref{3.24}, we consider two subcases:
\\
\textbf{\emph{Subcase i}}: $(\lambda_{N+1,N}V)(x)\neq0$. It then follows from \eqref{2.7} and \eqref{3.24} that for every $\rho=\rho(x)\in(0,1]$ we find  constants $u_{1}=u_{1}(x)$ and $u_{2}=-\rho u_{1}$, with $|u_{1}|$, $|u_{2}|$ being appropriately large, satisfying \eqref{3.7}.
\\
\textbf{\emph{Subcase ii}}: $(\lambda_{N+1,N}V)(x)=0$. Then, \eqref{3.24} becomes
\begin{equation}
    \label{3.25}
    \begin{split}
 &m^{(N+1)}(0)=(\rho+1)^{N+1}(f^{N+1}V)(x)+\sum_{k=1}^{j}u_{1}^{k}\Big( \Pi_{N+1,k}(\rho;x)
 \\
 &+\rho^{N}(\rho+1)(\lambda_{N+1,k}V)(x) \Big)+\sum_{k=\frac{j+3}{2}}^{\frac{N-3}{2}}u_{1}^{2k-1}\Big( \Pi_{N+1,2k-1}(\rho;x)
 \\
 &+\rho^{N}(\rho+1)(\lambda_{N+1,2k-1}V)(x) \Big)
 \\
 &+u_{1}^{N-2}\Big( \Pi_{N+1,N-2}(\rho;x)+\rho^{N}(\rho+1)(\lambda_{N+1,N-2}V)(x)\Big)
\end{split}
\end{equation}
\textbf{Step 2}: Since, according to \eqref{3.10b} with $n=N+1$ and $i=N-2$ the polynomials $\Pi_{N+1,N-2}(\rho;x)$ and $q(\rho)=\rho^{N}(\rho+1)$ are linearly independent, we distinguish two subcases:
\\
\textbf{\emph{Subcase i}}: There exists a constant $\rho=\rho(x)\in(0,1]$ such that $\Pi_{N+1,N-2}(\rho;x)+\rho^{N}(\rho+1)(\lambda_{N+1,N-2}V)(x)\neq0$. The latter in conjunction with \eqref{2.7} and \eqref{3.25} assert that there exist constants $u_{1}=u_{1}(x)$, $u_{2}=-\rho u_{1}$ with $|u_{1}|$, $|u_{2}|$ sufficiently large satisfying \eqref{3.7}.
\\
\textbf{\emph{Subcase ii}}: $\Pi_{N+1,N-2}(\rho;x)=0,\,\forall \rho\in\mathbb{R};(\lambda_{N+1,N-2}V)(x)=0$. Then \eqref{3.25} becomes:
\begin{equation}
    \label{3.26}
\begin{split}
 &m^{(N+1)}(0)=(\rho+1)^{N+1}(f^{N+1}V)(x)+\sum_{k=1}^{j}u_{1}^{k}\Big(
 \Pi_{N+1,k}(\rho;x)
 \\
 &+\rho^{N}(\rho+1)(\lambda_{N+1,k}V)(x) \Big)+\sum_{k=\frac{j+3}{2}}^{\frac{N-5}{2}}u_{1}^{2k-1}\Big( \Pi_{N+1,2k-1}(\rho;x)
 \\
 &+\rho^{N}(\rho+1)(\lambda_{N+1,2k-1}V)(x) \Big)
 \\
 &+u_{1}^{N-4}\Big( \Pi_{N+1,N-4}(\rho;x)+\rho^{N}(\rho+1)(\lambda_{N+1,N-4}V)(x)\Big)
 \\
\end{split}
\end{equation}
By taking into account \eqref{3.26} we may proceed quite similarly, as in Step 1 by induction and conclude that, either \eqref{3.7} holds for some $u_{1}=u_{1}(x)$, $u_{2}=-\rho u_{1}$ with $|u_{1}|$, $|u_{2}|$ sufficiently large, or
\begin{subequations}
\label{3.27}
\begin{equation}
    \label{3.27a}
    \begin{split}
    &\Pi_{N+1,i}(\rho;x)=0,\,\,\forall\,\rho\in\mathbb{R},\\
    &\text{for every odd}\,i\in\{j+2,...,N-2\}
    \end{split}
\end{equation}
\begin{equation}
    \label{3.27b}
    (\lambda_{N+1,i}V)(x)=0,\,\,\text{for every odd}\,i\in\{j+2,...,N\}
\end{equation}
\end{subequations}
and the procedure is terminated at the Step $(N-j+2)/2$ below, where, according to \eqref{3.27a} and \eqref{3.27b}, the original  expression \eqref{3.15} is written:
\begin{equation}
    \label{3.28}
\begin{split}
 &m^{(N+1)}(0)=(\rho+1)^{N+1}(f^{N+1}V)(x)
 \\
 &+\sum_{k=1}^{j-1}u_{1}^{k}\Big( \Pi_{N+1,k}(\rho;x)+\rho^{N}(\rho+1)(\lambda_{N+1,k}V)(x) \Big)
 \\
 &+u_{1}^{j}\left( \Pi_{N+1,j}(\rho;x)+\rho^{N}(\rho+1)(\lambda_{N+1,j}V)(x)\right)
 \end{split}
\end{equation}
\textbf{Step} \boldmath${(N-j+2)/2}$ \unboldmath: (notice that, since both $N$, $j$ are odd, $(N-j+2)/2$ is a positive integer): Only one case is considered in this last step, particularly, due to \eqref{2.8}, we have $(\lambda_{N+1,j}V)(x)\neq0$ and the latter in conjunction to \eqref{2.7}, together with validity of \eqref{3.10b} with $n=N+1$, $i=j$ and \eqref{3.28} imply that for every $\rho=\rho(x)\in(0,1]$, there exist constants $u_{1}=u_{1}(x)$, $u_{2}=-\rho u_{1}$ with $|u_{1}|$, $|u_{2}|$ sufficiently large such that inequality \eqref{3.7} holds.
\\
\textbf{CASE 4}: Suppose that conditions \eqref{2.5}, \eqref{2.7} and property P2(iii) hold. Namely, there exists an even integer $N\geq2$  satisfying \eqref{2.5}, \eqref{2.7} and \eqref{2.10}. It then follows from \eqref{3.9c} with $n=N+1$ that for every $\rho=\rho(x)\in(0,1]$, there exist constants $u_{1}=u_{1}(x)$, $u_{2}=-\rho u_{1}$ with $|u_{1}|$, $|u_{2}|$ sufficiently large such that inequality \eqref{3.7} is fulfilled.
\par By exploiting \eqref{3.1}-\eqref{3.7} and arguing as in the proof of [24, Proposition 3]  it follows that for every $\sigma>0$ and $x\neq0$ there exist a constant $\varepsilon=\varepsilon(x)\in(0,\sigma]$ and a measurable and locally essentially bounded control $u(\cdot,x):[0,\varepsilon]\to\mathbb{R}^{\text{m}}$ satisfying \eqref{1.3a}, \eqref{1.3b} with $a(s):=2s$ which asserts that \eqref{1.2} is SDF-SGAS. For completeness, we note that for the case where $(gV)(x)=0,\,\,(fV)(x)\leq0,\,\,x\neq0$, the corresponding control is defined as follows. Consider first the constants $\rho=\rho(x)\in(0,1]$, $u_{1}=u_{1}(x)$, $u_{2}=-\rho u_{1}$ as determined above (CASES 1-4). Let:
\begin{equation}
\label{3.29}
    \eta(s;t,x):=\begin{cases}
u_{2}=-\rho u_{1}, & s\in[0,t] \\ 
u_{1}, &s\in(t,t+\rho t] 
\end{cases}
\end{equation}
and for  every sufficiently small $\sigma=\sigma(x)>0$ and $\varepsilon\in(0,\sigma]$ define: 
\begin{equation}
    \label{3.30}
    u(\cdot,x):=\eta(\cdot;(\varepsilon/(1+\rho)),x)
\end{equation}
It then follows by taking into account \eqref{3.3}-\eqref{3.7}, \eqref{3.29}, \eqref{3.30} that the corresponding trajectory $\pi(\cdot,0,x,u)$ of \eqref{1.2} satisfies $V(\pi(\varepsilon,0,x,u(\cdot,x)))<V(x)$ and simultaneously $V(\pi(s,0,x,u(\cdot,x)))\leq 2V(x)$, for every $s\in(0,\varepsilon]$, which in conjunction with \eqref{3.1} and \eqref{3.2} imply both \eqref{1.3a}, \eqref{1.3b}.
\\
\textbf{Proof of Proposition 3}: Let $\Omega\subset\mathbb{R}^{\text{n}}$  be a nonempty bounded neighborhood of $0\in\mathbb{R}^{\text{n}}$. By \eqref{3.1} for the candidate $u$ it follows from assumption \eqref{2.11} that there exist a pair of continuous nonnegative functions $\theta:\mathbb{R}^{\text{n}}\to\mathbb{R}^{+}$ and $\xi:\mathbb{R}^{+}\to\mathbb{R}^{+}$ and a constant $C>0$ such that
        $|u(x)|\leq\left (\frac{|(fV)(x)+\theta(x)|}{|(gV)(x)|^{2}}+1 \right)|(gV)(x)|\leq\xi(|x|)+|(gV)(x)|\leq C,\,\,\forall x\in\Omega:\,\text{either}\,(gV)(x)\neq0,\,\text{or}\,(gV)(x)=0,\,\text{and}\,(fV)(x)<0$
 We next examine three cases:
 \\
\textbf{CASE 1}: Consider those nonzero $x$ for which $(gV)(x)=0$ and assumptions \eqref{2.5}, \eqref{2.7} hold along with one of the properties P1, P2(i). Then, by applying the same procedure with this used in the proof of Proposition 2, we conclude that, for every $\mu>0$ a vector $u=(u_{1},u_{2})\in\mathbb{R}^{2}$ can be determined such that \eqref{3.6} and \eqref{3.7} hold and further  
\begin{subequations}
\label{3.32}
\begin{equation}
    \label{3.32a}
    |u(x)|\leq\mu
\end{equation}
\begin{equation}
    \label{3.32b}
    \begin{split}
        &\forall x\in\Omega\setminus\{0\}:(gV)(x)=0,
        \\
        &\text{\emph{provided that either P1 or P2(i) holds}}
    \end{split}
\end{equation}
\end{subequations}
We next show that, for the remaining cases of the statement of Proposition 3, the same properties are fulfilled, namely, for every $\mu>0$ and $\rho=\rho(x)\in(0,1]$ a vector  $u=(u_{1},u_{2})\in\mathbb{R}^{2}$ can be determined such that \eqref{3.6}, \eqref{3.7} and \eqref{3.32a} hold.
\\
\textbf{CASE 2}: Consider those nonzero $x$ for which $(gV)(x)=0$ and assume that \eqref{2.5}, \eqref{2.7} and Property P2(ii) are fulfilled with $j=j(N)=N$ for some odd integer $N=N(x)\geq1$ and in such a way that the additional property \eqref{2.12} is fulfilled. Then, we may apply the same procedure with this used in the proof of Proposition 2 and prove \eqref{3.6}. We next establish that \eqref{3.7} holds as well for certain arbitrarily small $u_{1}=u_{1}(x)$, $u_{2}=-\rho u_{1}$ and $\rho\in(0,1]$. Notice that condition \eqref{2.12} in conjunction with \eqref{3.11} and \eqref{3.17} imply:
\begin{subequations}
\label{3.33}
\begin{equation}
    \label{3.33a}
    \Pi_{N+1,i}(\rho;x)=0,\,\,\text{for every}\,\,i\in\{1,...,N-1\}
\end{equation}
\begin{equation}
    \label{3.33b}
    (\lambda_{N+1,i}V)(x)=0,\,\,\text{for every}\,\,i\in\{1,...,N-1\}
\end{equation}
\end{subequations}
By taking into account \eqref{3.9c} with $n=N+1$ and \eqref{3.33a}, \eqref{3.33b} we have
\begin{equation}
    \label{3.34}
    \begin{split}
    m^{(N+1)}(0)&=(\rho+1)^{N+1}(f^{N+1}V)(x)
    \\
    &+u_{1}^{N}\rho^{N}(\rho+1)(\lambda_{N+1,N}V)(x)
    \end{split}
\end{equation}
It follows from \eqref{2.7}, \eqref{2.8}, \eqref{3.34} and oddness of $N$ that for every $\mu>0$ and $\rho=\rho(x)\in(0,1]$ a vector $u=(u_{1},u_{2})\in\mathbb{R}^{2}$ can be found with $u_{1}=u_{1}(x)$, $u_{2}=-\rho u_{1}$ in such a way that \eqref{3.6}, \eqref{3.7} and \eqref{3.32a} are  fulfilled.
\\
\textbf{CASE 3}: Suppose finally that, in addition to \eqref{2.5}, \eqref{2.7}, both P2(iii) and condition \eqref{2.12} hold for some even integer $N=N(x)\geq2$. We again may apply the same procedure with this used in the proof of Proposition 2 and prove \eqref{3.6}. In order to show \eqref{3.7}, we again notice that, due to \eqref{2.12}, \eqref{3.11} and \eqref{3.17}, the equalities \eqref{3.33a}, \eqref{3.33b} hold with $N=N(x)\geq2$ even. Consequently, by taking into account \eqref{3.9c} with $n=N+1$ and \eqref{3.33a}, \eqref{3.33b}, we have:
\begin{equation}
\label{3.35}
\begin{split}
     m^{(N+1)}(0)&=(\rho+1)^{N+1}(f^{N+1}V)(x)
     \\
     &+u_{1}^{N}\rho^{N}(\rho+1)(\lambda_{N+1,N}V)(x),\,\,N\geq2
     \end{split}
\end{equation}
(where now $N\geq2$ is even). We conclude, by taking into account \eqref{2.7}, \eqref{2.10} and \eqref{3.35}, that for every $\mu>0$ and $\rho=\rho(x)\in(0,1]$ a vector $u=(u_{1},u_{2})\in\mathbb{R}^{2}$ can be found with $u_{1}=u_{1}(x)$, $u_{2}=-\rho u_{1}$ in such a way that \eqref{3.6}, \eqref{3.7} and \eqref{3.32a} are  fulfilled.
\par By invoking \eqref{3.1}-\eqref{3.7}, the conclusions of Cases 1-3 above and \eqref{3.32} we may establish as in  the proof of Proposition 2, that in all cases the desired \eqref{1.3a}, \eqref{1.3b} together with \eqref{1.4} are fulfilled, hence, according to Proposition 1,  system \eqref{1.2} is BSDF-SGAS.
\section{An illustrative example}
The next example illustrates the nature of the results of the previous section. Consider the system
\begin{subequations}
\label{4.1}
\begin{equation}
\label{4.1a}
    \begin{pmatrix}
\dot{x}\\ 
\dot{y}
\end{pmatrix}=f+ug,\,\,x\in\mathbb{R}^{\text{n}},\,\,y\in\mathbb{R},\,\,u\in\mathbb{R}
\end{equation}
\text{where the mappings $f,g:\mathbb{R}^{\text{n}+1}\to\mathbb{R}^{\text{n}+1}$ have the form}
\begin{equation}
    \label{4.1b}
    \begin{split}
    f(x,y)&=\begin{pmatrix}
a(x)+y\beta(x)+y^{2}\gamma(x)+y^{3}\delta(x)\\ 
0
\end{pmatrix},\\
g(x,y)&=\begin{pmatrix}
0\\ 
1
\end{pmatrix},\,x\in\mathbb{R}^{\text{n}},\,\,y\in\mathbb{R}
\end{split}
\end{equation}
\end{subequations}
and the functions $a:\mathbb{R}^{\text{n}}\to\mathbb{R}$, $\beta:\mathbb{R}^{\text{n}}\to\mathbb{R}$,  $\gamma:\mathbb{R}^{\text{n}}\to\mathbb{R}$, $\delta:\mathbb{R}^{\text{n}}\to\mathbb{R}$ are smooth, with $a(0)=0$. We assume that there exists a smooth function $W:\mathbb{R}^{\text{n}}\to\mathbb{R}^{+}$ being positive definite and proper such that, if we define the sets:
\begin{subequations}
\label{4.2}
\begin{equation}
\label{4.2a}
    E_{1}:=\{x\in\mathbb{R}^{\text{n}}\setminus\{0\}:(aW)(x)<0\}
\end{equation}
\begin{equation}
    \label{4.2b}
    E_{2}:=\{x\in\mathbb{R}^{\text{n}}\setminus\{0\}:(aW)(x)\leq0,\,(\beta W)(x)\neq0\}
\end{equation}
\begin{equation}
    \label{4.2c}
    E_{3}:=\begin{Bmatrix}x\in\mathbb{R}^{\text{n}}\setminus\{0\}:
    (aW)(x)=(\beta W)(x)=0,(\gamma W)(x)<0\end{Bmatrix}
\end{equation}
\begin{equation}
    \label{4.2d}
    E_{4}:=\begin{Bmatrix}x\in\mathbb{R}^{\text{n}}\setminus\{0\}:(aW)(x)=0,\\
    (\beta W)(x)=(\gamma W)(x)=0,
    (\delta W)(x)\neq0\end{Bmatrix}
\end{equation}
\begin{equation}
    \label{4.2e}
    \begin{split}
    E_{5}:=\begin{Bmatrix}
x\in\mathbb{R}^{\text{n}}\setminus\{0\}:(aW)(x)=(\beta W)(x)=0,\\
(\gamma W)(x)=(\delta W)(x)=([[a,\gamma],a]W)(x)=0,\\([a,\delta]W)(x)\neq0
\end{Bmatrix}
\end{split}
\end{equation}
\text{then, the following holds}
\begin{equation}
    \label{4.2f}
    E_{1}\cup E_{2}\cup E_{3}\cup E_{4}\cup E_{5}=\mathbb{R}^{\text{n}}\setminus\{0\}
\end{equation}
\end{subequations}
We also assume that each $E_{i}$, $i=1,2,3,4,5$ is nonempty and satisfies:
\begin{equation}
    \label{4.3}
    \begin{split}
        &0\in cl E_{i},\,i=1,2,3,4,5,\,(E_{1}\cup E_{2})\cap(E_{3}\cup E_{4}\cup E_{5})=\partial E_{1},
        \\
        &int E_{3}\neq\varnothing,\,cl(E_{4}\cup E_{5})\subset \partial E_{3},\,int(E_{4}\cup E_{5})=\varnothing
    \end{split}
\end{equation}
(see Fig. 1).
\begin{figure}[!t]
\centerline{\includegraphics[width=\columnwidth]{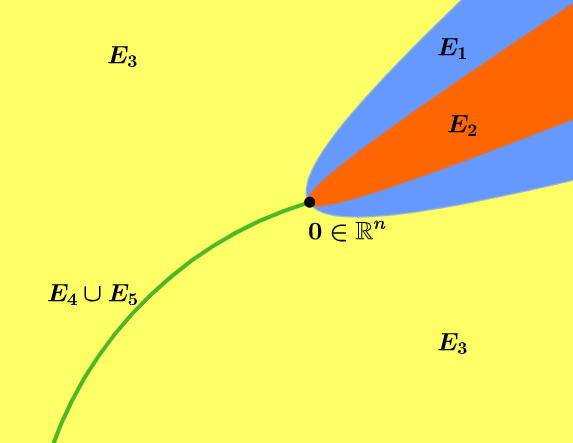}}
\caption{ Graphical representation of the sets $E_{i}$, $i=1,2,3,4,5$ }
\end{figure}
\\
CLAIM:  Under previous assumptions, the system \eqref{4.1} satisfies  the hypotheses of Propositions 2, 3 with
\begin{equation}
    \label{4.4}
    V(x,y):=W(x)+y^{2},\,\,x\in\mathbb{R}^{\text{n}},\,\,y\in\mathbb{R}
\end{equation}
 therefore, it is BSDF-SGAS.
 \\
\emph{Proof of Claim} : Obviously, $V:\mathbb{R}^{\text{n}}\times\mathbb{R}\to\mathbb{R}^{+}$, as defined by \eqref{4.4} is smooth, positive definite and proper. Also, according to definitions \eqref{4.1b} and \eqref{4.4}, it holds that
\begin{equation}
\label{4.5}
\begin{split}
    (fV)(x,y)&:=(aW)(x)+y(\beta W)(x)+y^{2}(\gamma W)(x)
    \\
    &+y^{3}(\delta W)(x),\, \forall (x,y)\in\mathbb{R}^{\text{n}}\times\mathbb{R}
    \end{split}
\end{equation}  
and due to \eqref{4.2}, the map $aW$ takes nonpositive values. It follows from \eqref{4.4} and \eqref{4.5}, that condition \eqref{2.11} is satisfied with $(x,y)$ instead of $x$, $V$ as given by \eqref{4.4}, $\theta(x,y)=-(aW)(x)$, $(gV)(x):=2y$ and for certain $\xi:\mathbb{R}^{+}\to\mathbb{R}^{+}$ in such a way that
$|y(\beta W)(x)+y^{2}(\gamma W)(x)+y^{3}(\delta W)(x)|\leq |y|\xi(|(x,y)|),\,\forall (x,y)\in\mathbb{R}^{\text{n}}\times\mathbb{R}$
Assume next that for some nonzero $(x,y)$ we have:
\begin{equation}
    \label{4.6}
    (gV)(x,y)=0\overset{(4.1\text{b}),(4.4)}{\Leftrightarrow }x\neq0,\,y=0
\end{equation}
According to \eqref{4.2} and \eqref{4.3}, we distinguish five cases concerning \eqref{4.6}:
\\
\textbf{CASE 1}: $x\in E_{1}$, $y=0$. Then condition \eqref{2.4} of Proposition 2 is fulfilled for the dynamics $f,g$ of system \eqref{4.1}, i.e.:
\begin{equation}
    \label{4.7}
    (gV)(x,y)=0\Rightarrow (fV)(x,y)<0,\,\text{for}\,(x,y)\in E_{1}\times\{0\}
\end{equation}
Indeed, from \eqref{4.1b} and \eqref{4.4} we find $(fV)(x,y)=(DVf)(x,y)=(aW)(x)+y(\beta W)(x)+y^{2}(\gamma W)(x)+y^{3}(\delta W)(x)$ and recalling \eqref{4.6} it follows that $(fV)(x,y)=(aW)(x)$ for every $(x,y)\in\mathbb{R}^{\text{n}}\times\{0\}$. Since $x\in E_{1}$, the latter implies \eqref{4.7}.
\\
\textbf{CASE 2(i)}: $x\in E_{2}$, $y=0$; $(aW)(x)=0$. Under previous assumption, we claim that condition \eqref{2.5} together with Property P2(i) are satisfied. Particularly, the following hold: Condition \eqref{2.5} with $N=1$, i.e.: $(fV)(x,y)=0$. Equality in condition \eqref{2.7} with $N=1$, particularly: $(f^{2}V)(x,y)=0$. Condition \eqref{2.8} is fulfilled with $N=1$ and (odd) $j=1$, i.e. $(\lambda_{2,1}V)(x,y)\neq0$.  Indeed, $aW$ is smooth and according to \eqref{4.2}, in our case, it takes maximum value. Consequently, we have $(aW)(x)=0$, $D(aW)(x)=0$, thus, by \eqref{4.4}, $(f^{i}V)(x,y)=(a^{i}W)(x)$, $i=1,2,...$, therefore, $(f^{i}V)(x,y)=0$, $i=1,2,...$, $x\in E_{2}$, $y=0$.  Also, by taking into account definition \eqref{4.4} and \eqref{4.6} we find $(\lambda_{2,1}V)(x,y)=-(\beta W)(x)\neq0$, for $x\in E_{2}$, $y=0$.
\\
\textbf{CASE 2(ii)}: $x\in E_{2}$, $y=0$; $(aW)(x)<0$. Under previous assumption, we show as in Case I that \eqref{4.7} holds.
\\
\textbf{CASE 3}: $x\in E_{3}$, $y=0$. We claim that in this case, condition \eqref{2.5}, together with Property P2(iii) and \eqref{2.12} are fulfilled; particularly, the following hold: condition \eqref{2.5} with $N=2$, i.e., $(f^{i}V)(x,y)=0$, $i=1,2$; $(\lambda_{2,1}V)(x,y)=0$, equality in condition \eqref{2.7} with $N=2$, i.e., $(f^{3}V)(x,y)=0$, inequality \eqref{2.10} with $N=2$, i.e. $(\lambda_{3,2}V)(x,y)<0$ and further \eqref{2.12} is fulfilled with $N=2$, i.e., $(\lambda_{3,1}V)(x,y)=(\lambda_{2,1}\lambda_{1,0}V)(x,y)=(\lambda_{1,0}\lambda_{2,1}V)(x,y)=0$. Indeed, by first recalling \eqref{4.2c}-\eqref{4.2f}, rest properties \eqref{4.3}  and smoothness of $aW$, we find $(aW)(x)=0$, $\forall x\in\partial E_{1}\cup cl(E_{3}\cup E_{4}\cup E_{5})$. It follows that for $x\in E_{3}$, $y=0$ we have $(f^{i}V)(x,y)=(a^{i}W)(x)=0$, $i=1,2,...$. Moreover, it holds $(\beta W)(x)=0$, $\forall x\in\partial E_{1}\cup cl(E_{3}\cup E_{4}\cup E_{5})$. It follows that $(\lambda_{2,1}V)(x,y)=0$, for $x\in E_{3}$, $y=0$. The above asserts that \eqref{2.5} and equality in \eqref{2.7} hold with $N=2$. Also, by taking into account \eqref{2.1} and \eqref{4.6} we have $(\lambda_{3,2}V)(x)=2(\gamma W)(x)$ and  therefore from \eqref{4.2c} we get $(\lambda_{3,2}V)(x,y)<0$, i.e., inequality \eqref{2.10} of Proposition 2 holds with $N=2$. Finally, by again recalling \eqref{4.2} and \eqref{4.3} we have $(a\beta W)(x)=(\beta aW)(x)=([a,\beta]W)(x)=0$, $\forall x\in E_{3}$, which implies  validity of \eqref{2.12} (Proposition 3) with $N=2$.
 \\
 \textbf{CASE 4}: $x\in E_{4}$, $y=0$. We claim that, when $x\in E_{4}$, $y=0$, then condition \eqref{2.5} together with Property P2(ii) and \eqref{2.9a} of Proposition 2 are satisfied. Particularly, condition \eqref{2.5} holds with $N=3$, i.e., $(f^{i}V)(x,y)=0$, $i=1,2,3$, equality in condition \eqref{2.7} holds with $N=3$, i.e., $(f^{4}V)(x,y)=0$,  property \eqref{2.8} is fulfilled with $N=3$ and (odd) $j=3$, i.e., $(\lambda_{4,3}V)(x,y)\neq0$ and finally, condition \eqref{2.9a} is satisfied with $N=3$, (odd) $j=3$ and (even) $q=2$.
Indeed, by recalling \eqref{4.6}, we get $(f^{i}V)(x,y)=(a^{i}V)(x)=0$, $i=1,2,3,4$, $(a\beta W)(x)=(\beta aW)(x)=0$, $\forall x\in E_{4}$. It follows that for the case $x\in E_{4}$, $y=0$ it holds $(\lambda_{2,1}V)(x,y)=(\lambda_{3,1}V)(x,y)=(\lambda_{1,0}\lambda_{2,1}V)(x,y)=(\lambda_{2,1}\lambda_{1,0}V)(x,y)=0$. Moreover, according to \eqref{4.2d}, we have $(\lambda_{4,3}V)(x,y)=-6(\delta W)(x)\neq0$ and the function $\gamma W$ takes maximum on the region $E_{4}\cup E_{5}$, hence, $(\gamma W)(x)=0$, $D(\gamma W)(x)=0$, for $x\in E_{4}$, $y=0$  and this implies $(\lambda_{3,2}V)(x,y)=2(\gamma W)(x)=(a\gamma W)(x)=(\gamma a W)(x)=0$, for $x\in E_{4}$, $y=0$. By taking into account the previous facts we can easily verify that \eqref{2.5} is satisfied with $N=3$, as well as \eqref{2.8} and \eqref{2.9a}  with $N=3$, $j=3$ and (even) $q=2$ respectively.
\\
\textbf{CASE 5}: $x\in E_{5}$, $y=0$. We claim that \eqref{2.5} together with Property P2(i) are satisfied. Particularly, implication \eqref{2.5} holds with $N=4$.
Also, equality in \eqref{2.7} with $N=4$ is satisfied, i.e., $(f^{5}W)(x,y)=0$, property \eqref{2.8} holds with $N=4$ and $j=3$, i.e., $(\lambda_{5,3}V)(x,y)\neq0$ and finally \eqref{2.9a}, \eqref{2.9b} hold with $N=4$, $j=3$ and $q=2$.
The establishment of the above facts is a consequence of the fact that $(a^{i}W)(x)=0$, $i=1,2,...$, $([\beta,\gamma]W)(x)=0$, $(\lambda_{5,3}V)(x,y)=6([a,\delta]W)(x)\neq0$.
Details are left to the reader. 

\section{conclusion}
This paper presents new results on sampled-data feedback stabilization for affine in the control of nonlinear systems with nonzero drift term under the presence of a control Lyapunov function associated with Lie algebraic hypotheses concerning the dynamics of the system.


\begin{thebibliography}{00}
\bibitem{art:1}{F. Ancona and A. Bressan, "Patchy vector fields and asymptotic stabilization," \emph{ESAIM-COCV} vol.4, pp.445-471, 1999.}

\bibitem{art:2}{Z. Artstein, "Stabilization with relaxed controls," \emph{Nonlinear Analysis TMA}, vol.7, pp. 1163-1173, 1983.}

\bibitem{art:3}{A. Bacciotti and L. Mazzi, "From Artstein-Sontag Theorem to the min-projection strategy," \emph{Trans. of the Institute of Measurement and Control}, vol.32, no.6, pp. 571-581, 2010.}

\bibitem{art:4}{F.H. Clarke, Y.S. Ledyaev, E.D. Sontag and A.I. Subbotin, "Asymptotic controllability implies feedback stabilization," \emph{IEEE Trans. Autom. Control}, vol. 42, no. 10, pp. 1394-1407, 1997.}

\bibitem{art:5}{F.H. Clarke, Y.S. Ledyaev, L. Rifford and R.J. Stern, "Feedback stabilization and Lyapunov functions," \emph{SIAM J.  Control Optim.}, vol. 39, no. 1,pp. 25-48, 2000.}

\bibitem{art:6}{R. Goebel and A.R. Teel, "Direct design of robustly asymptotically stabilizing hybrid feedback," \emph{ESAIM-COCV}, vol. 15, no. 1, pp. 205-213, 2009.}

\bibitem{art:7}{L. Gr\"{u}ne and D. Ne\v{s}i\'{c}, "Optimization based stabilization of sampled-data nonlinear systems via their approximate discrete-time models," \emph{SIAM J. Control  Optim.}, vol. 42, pp. 98-122, 2003. }

\bibitem{art:8}{H.G. Hermes , "Controlled Stability," \emph{Ann. Mat. Pura. Appl.}, IV144, pp. 103-119, 1977.}

\bibitem{art:9}{I. Karafyllis, "Stabilization by Means of Time-Varying Hybrid Feedback," \emph{Mathematics of Control, Signals and Systems}, vol. 18, no. 3, pp. 236-259, 2006.}

\bibitem{art:10}{W. Lin and W. Wei, "Semi-Global Asymptotic Control by Sampled-Data Output Feedback," \emph{IFAC}, vol. 51, no.18, pp. 596-601, 2018.}

\bibitem{art:11}{N. Marchand and M. Alamir, \emph{Asymptotic controllability implies continuous discrete-time feedback stabilization}.  Nonlinear Control in the Year 2000, vol. 2, Springer, Berlin, Heidelberg, New York, 2000.}

\bibitem{art:12}{S. Monaco, D. Normand-Cyrot and M. Mattioni, "Sampled-Data Stabilization of Nonlinear Dynamics With Input Delays Through Immersion and Invariance," \emph{IEEE Transactions on Automatic Control}, vol. 62, no. 5, pp. 1-1, 2016.}

\bibitem{art:13}{M. Motta and F. Rampazzo, "Asymptotic Controllability and Lyapunov-like Functions Determined by Lie Brackets," \emph{SIAM Journal on Control and Optimization}, vol. 56, no. 2, pp. 1508-1534, 2018.}

\bibitem{art:14}{D. Ne\v{s}i\'{c}, A.R. Teel and P.V. Kokotovic, "Sufficient conditions for stabilization of sampled-data nonlinear systems via discrete-time approximations," \emph{Systems and Control Lett.}, vol. 38, no. 4-5, pp. 259-270, 1999.}

\bibitem{art:15}{D. Ne\v{s}i\'{c} and A.R. Teel, "A framework for stabilization of nonlinear sampled-data systems based on their approximate discrete-time models," \emph{IEEE Trans. Autom. Control}, vol. 49, pp. 1103-1122, 2004.}

\bibitem{art:16}{C. Prieur, "Asymptotic controllability and robust asymptotic stabilizability," \emph{SIAM J. Control Optim.}, vol.43, pp. 1888-1912, 2005.}

\bibitem{art:17}{C. Prieur, R. Goebel and A.R. Teel, "Hybrid feedback control and robust stabilization of nonlinear systems," \emph{IEEE Transactions on Automatic Control}, vol. 52, no. 11, pp. 2103 - 2117, 2007.}

\bibitem{art:18}{H. Shim, A.R. Teel, "Asymptotic controllability and observability imply semiglobal practical asymptotic stabilizability by sampled-data output feedback," \emph{Automatica}, vol. 39, pp. 441-454, 2003.}

\bibitem{art:19}{E.D. Sontag, "A "universal" construction of Artstein's theorem on nonlinear stabilization,"   \emph{Systems and Control Lett.}.  vol. 13, pp. 117-123, 1989.}


\bibitem{art:20}{J. Tsinias, “Sufficient Lyapunov-like conditions for stabilization,” \emph{Math. Contr. Sign. Syst.}, vol. 2, pp. 343-357, 1989.}

\bibitem{art:21}{J. Tsinias, “Remarks on asymptotic controllability and sampled-data feedback stabilization for autonomous systems,” \emph{IEEE Trans. Autom. Control}, vol. 55, pp. 721-726, 2010.}

\bibitem{art:22}{J. Tsinias, “Small-gain type sufficient conditions for sampled-data feedback stabilization for autonomous composite systems,” \emph{IEEE Trans. Autom. Control}, vol. 56, pp. 1725-1729, 2011.}

\bibitem{art:23}{J. Tsinias, “New results on sampled-data feedback  stabilization for autonomous nonlinear systems,” \emph{Systems and Control Lett.}, vol. 61, pp. 1032-1040, 2012.}

\bibitem{art:24}{J. Tsinias and D.Theodosis, "Sufficient Lie Algebraic Conditions for Sampled-Data Feedback Stabilizability of Affine in the Control Nonlinear Systems," \emph{IEEE Transactions on Automatic Control}, vol. 61, pp. 1334-1339, 2016.}
\end{thebibliography}
\end{document}